\documentclass[11pt]{article}
\usepackage{graphicx}
\usepackage{amsmath}
\usepackage{amssymb}
\usepackage{amsfonts}
\usepackage{mathtools}
\usepackage[numbers]{natbib}
\usepackage{url}

\newtheorem{theorem}{Theorem}

\def\R{\mathbb{R}}

\def\xit{{\xi_*}}

\title{Continuation and bifurcations of periodic orbits and symbolic dynamics in the Swift-Hohenberg equation}
\author{
    Jakub Czwórnóg and Daniel Wilczak\\
    Faculty of Mathematics and Computer Science\\
    Jagiellonian Univeristy\\
    Łojasiewicza 6, 30--348 Kraków, Poland\\
    \texttt{\{kuba.czwornog,Daniel.Wilczak\}@uj.edu.pl}
}
\date{\today}

\begin{document}
\maketitle
\begin{abstract}
  Steady states of the Swift--Hohenberg \cite{SwiftHohenberg} equation are studied. For the associated four--dimensional ODE we prove that on the energy level $E=0$ two smooth branches of even periodic solutions are created through the saddle-node bifurcation. We also show that these orbits satisfy certain geometric properties, which implies that the system has positive topological entropy for an explicit and wide range of parameter values of the system. 
  
  The proof is computer-assisted and it uses rigorous computation of bounds on certain Poincar\'e map and its higher order derivatives.
\end{abstract}

\section{Introduction}
The Swift-Hohenberg equation \cite{SwiftHohenberg} is a fundamental partial differential equation that plays a crucial role in the study of pattern formation and models various phenomena in physics and biology \cite{TLIDI19941475,PhysRevLett.73.2978,MERON201270}. Originally the one-dimensional equation has the form
\begin{equation} \label{eq:SH_equation}
	\frac{\partial U}{\partial t} = 
	- \left ( \frac{\partial^2}{\partial X^2} + 1 \right )^2 U + \alpha U - U^3.
\end{equation}
The existence of some types of stationary solutions such as periodic and homoclinic  \cite{Glebsky,BurkeKnoblach,Deng} in one-dimensional case is studied analytically. There are also many numerical simulations that demonstrate the existence of complicated structures and behaviour \cite{Yang,Su,sanchez2013numerical}. However, many observed phenomena are not proved by means of mathematical rigour. 

The aim of this article is to reproduce and extend results from \cite{chaos} about one-dimensional stationary Swift-Hohenberg \cite{SwiftHohenberg} equation

\begin{equation} \label{eq:SH_ODE}
	-U'''' - 2 U'' + (\alpha - 1)U - U^3 = 0.
\end{equation}
This equation has conserved energy
\begin{equation}
	E = U''' U' - \frac{1}{2} (U'')^2 + U'^2 - \frac{\alpha - 1}{2} U^2 + \frac{1}{4} U^4 + \frac{(\alpha - 1)^2}{4}.
\end{equation}
The authors of \cite{chaos} proved the following theorem.
\begin{theorem} \cite[Proposition 1]{chaos} \label{thm:vandenBergLessard}
	The dynamics of the ODE \eqref{eq:SH_ODE} on the energy level $E=0$ is chaotic for all $\alpha \geq 2$ in the sense that certain Poincar\'e map is semi-conjugated to a subshift of finite type with positive topological entropy.
\end{theorem}
The proof in \cite{chaos} splits into two parts. First, the authors prove \cite[Theorem 2]{chaos} that the existence of a periodic orbits of (\ref{eq:SH_ODE}) with the parameter value $\alpha>\frac{3}{2}$ on the energy level $E=0$ and satisfying certain geometric properties (see Section~\ref{sec:equationForPOs} for details) implies the existence of symbolic dynamics for (\ref{eq:SH_ODE}). Then in \cite[Theorem 3]{chaos} the authors show, that the assumptions of \cite[Theorem 2]{chaos} are checked for all $\alpha\geq 2$. The proof of \cite[Theorem 3]{chaos} is computer-assisted and it is based on so-called radii polynomial approach. The authors reformulate the problem of the existence of a periodic orbit as a zero-finding problem for the Fourier coefficients of the orbit to be found. Then, by means of interval arithmetics \cite{Moore1966,IEEE1788-2015} and Newton-Kantorowitch type argument, it is shown that this infinite-dimensional system of algebraic equations has a branch of isolated solution parameterized by  $\alpha\geq 2$. 

The aim of this paper is to prove the following result.
\begin{theorem}\label{thm:main}
	There exists $\alpha^*\in 1.9690842080_{101989}^{293001}$ and two smooth curves $U_\pm:[\alpha^*,\infty)\to\mathbb R^4$ such that 
	\begin{enumerate}
		\item $E\left(U_\pm(\alpha)\right)=0$ for $\alpha\geq \alpha^*$, 
		\item the solution to \eqref{eq:SH_ODE} with the parameter value $\alpha$ and with the initial condition $(U,U',U'',U''')=U_\pm(\alpha)$ is an even and periodic function satisfying geometric properties from \cite[Theorem 2]{chaos},
		\item $U_-(\alpha)\neq U_+(\alpha)$ for $\alpha>\alpha^*$ and $U_-(\alpha^*)=U_+(\alpha^*)$,
		\item the fold bifurcation occurs at $\alpha^*$.
	\end{enumerate}
\end{theorem}
For the proof of Theorem~\ref{thm:main} we reformulate the problem of the existence of an even periodic orbit as a zero-finding problem for an univariate map (certain Poincar\'e map). In that sense, the proposed approach is much simpler and geometric as we work directly in the phase space of the ODE rather than an infinite dimensional space of Fourier coefficients. Moreover, computer-assisted verification of the result is significantly faster than the one presented in \cite{chaos} -- less than 13 minutes for two branches of periodic orbits and the fold bifurcation versus 12 hours for one branch.

The paper is organised as follows. In Section~\ref{sec:forcingTheorem} we recall the required geometric properties of these periodic orbits from \cite[Theorem 2]{chaos}. In Section~\ref{sec:equationForPOs} and derive a simple scalar equation for even periodic orbits of \eqref{eq:SH_ODE}. In Section~\ref{sec:cap} we present details of the computer-assisted proof of Theorem~\ref{thm:main}.

\section{Symbolic dynamics from a single periodic orbit}
\label{sec:forcingTheorem}
For self-consistency of the article we recall here \emph{the forcing theorem} \cite[Theorem 2]{chaos} and required geometric conditions for periodic orbits.

The system \eqref{eq:SH_ODE} has two equilibria $\pm\sqrt{\alpha-1}$ on the energy level $E=0$. For $1<\alpha\leq \frac{3}{2}$ they are stable foci (pure imaginary eigenvalues) and for all $\alpha>\frac{3}{2}$ they are of saddle-focus type. Thus, for $\alpha>\frac{3}{2}$ near these equilibria families of hyperbolic (on isolated energy level) periodic orbits appear giving rise to more complicated dynamics. 

In order to make the interval of parameters bounded, and thus much easier for computer analysis, following the ideas from \cite{chaos} we perform a change of coordinates 
\begin{equation}\label{eq:changeVariables}
	y = \frac{X}{\sqrt[4]{\alpha - 1}}, \quad u(y) = \frac{U(X)}{\sqrt{\alpha - 1}}, \quad \xi = \frac{2}{\sqrt{\alpha - 1}}.
\end{equation}
Now the parameter range $\alpha \geq \frac{3}{2}$ corresponds to $0<\xi \leq \sqrt8$ and \eqref{eq:SH_ODE} becomes 
\begin{equation} \label{eq:newSH_ODE}
	-u'''' - \xi u'' + u - u^3 = 0
\end{equation}
with the energy
\begin{equation}\label{eq:newSH_energy}
	E = u''' u' - \frac{1}{2}(u'')^2+ \frac{\xi}{2}(u')^2 + \frac{1}{4}(u^2 - 1)^2.
\end{equation}

\begin{figure}[htbp] 
	\centering
	\includegraphics[width=1\textwidth]{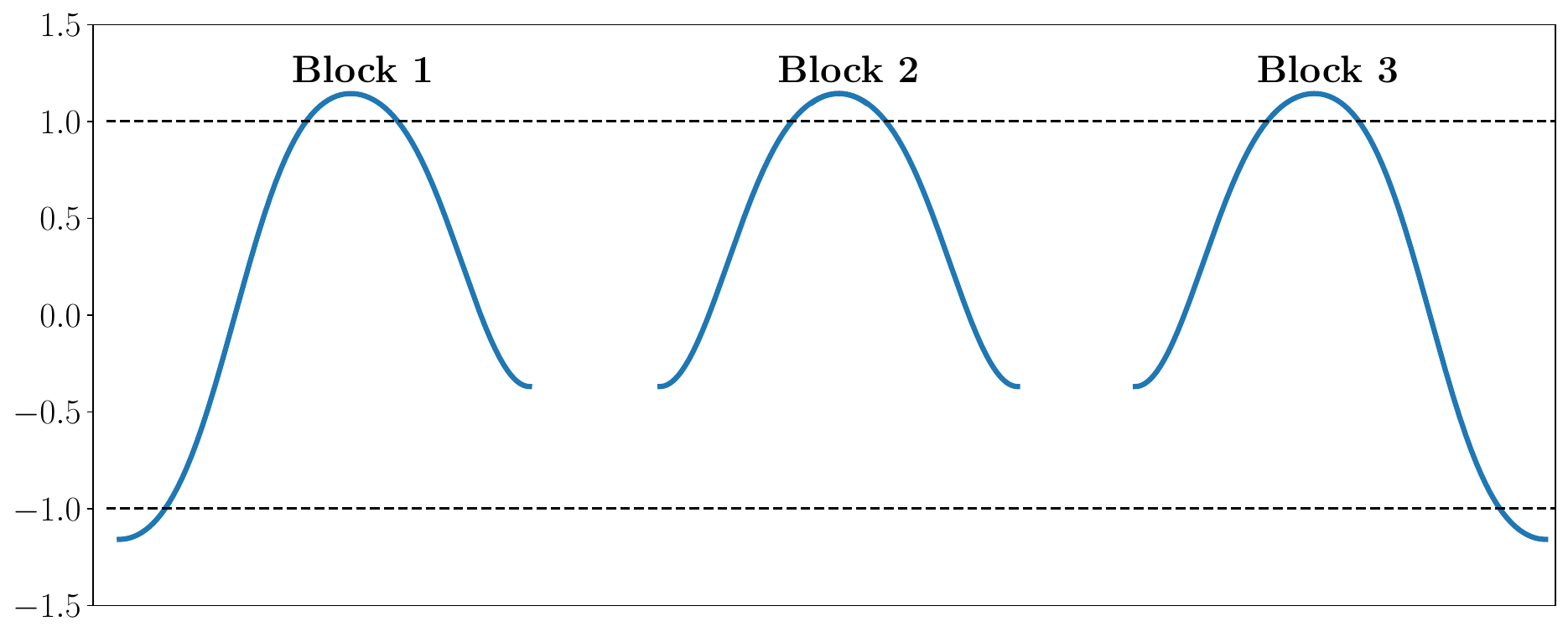}
	\caption{Three segments of a periodic orbit that lead to construction of symbolic dynamics. Blocks 1 and 2 can be followed by either Block 2 or Block 3. However, Block 3 must be followed by Block 1.}
	\label{fig:orbit}
\end{figure}

\begin{theorem} \cite[Theorem 2]{chaos} \label{theorem:forcing}
	Let $\xi \in [0,\sqrt{8})$ and suppose that there exists a periodic solution $\tilde{u}$ of \eqref{eq:newSH_ODE} at the energy level $E=0$, satisfying the following geometric conditions
	\begin{itemize}
		\item $\tilde{u}$ has exactly four monotone laps in one period and extrema $\tilde{u}_1,\tilde{u}_2,\tilde{u}_3,\tilde{u}_4$;
		\item $\tilde{u}_1, \tilde{u}_3$ are minima, and $\tilde{u}_2,\tilde{u}_4$ are maxima;
		\item $\tilde{u}_1 < -1 < \tilde{u}_3 < 1 < \tilde{u}_2, \tilde{u}_4$;
		\item $\tilde u$ is symmetric at its minima.
	\end{itemize}
	Then the system is chaotic in the sense that there exists a two-dimensional Poincar\'e return
	map which has a compact invariant set on which the topological entropy is positive.
\end{theorem}
\begin{figure}[htbp] 
	\centering
	\includegraphics[width=.48\textwidth]{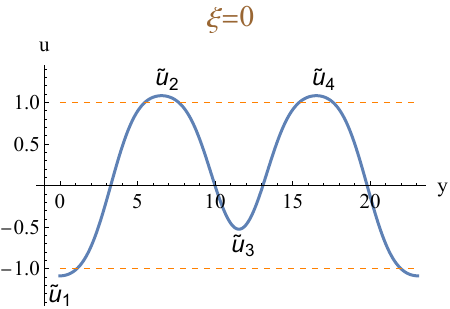}
	\includegraphics[width=.48\textwidth]{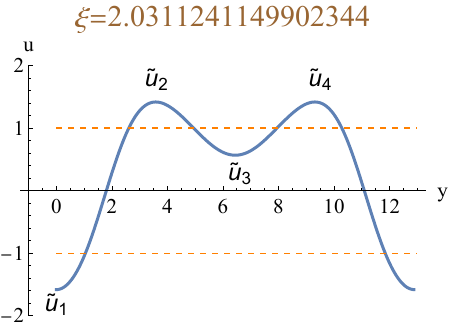}
	\caption{Shape of periodic orbits satisfying the geometric properties of Theorem~\ref{theorem:forcing}.}
	\label{fig:periodic_orbits}
\end{figure}
Geometric conditions about periodic orbits are illustrated in Fig.\ref{fig:periodic_orbits}. The solutions that belong to the chaotic invariant set resulting from Theorem~\ref{theorem:forcing} are encoded by their extrema. While the maximum is always bigger than $1$ we are free to choose minimum to be either less than $-1$ or between $-1$ and $1$. This leads to symbolic dynamics on three symbols, where each symbol is a building block of solution -- as illustrated in Fig.~\ref{fig:orbit}.

\section{Even periodic orbits as a zero-finding problem.}
\label{sec:equationForPOs}
Equation \eqref{eq:newSH_ODE} can rewritten as a system of first order ODEs
\begin{equation} \label{eq:SH_system}
	\begin{cases}
		x'(t) &= y(t) \\
		y'(t) &= z(t) \\
		z'(t) &= w(t) \\
		w'(t) &= -\xi z(t) + x(t) - x^3(t).
	\end{cases}
\end{equation}
Define a Poincar\'e section
\begin{equation*}
	\Pi = \left \{(x,y,z,w)  :  y = 0 \right \}.    
\end{equation*}
and let $\mathcal P_\xi:\Pi\to \Pi$ be the associated Poincar\'e map for \eqref{eq:SH_system} with the fixed parameter value $\xi$. The choice $y=u'=0$ of Poincar\'e  section is very natural when we are looking for periodic orbits with fixed number of extrema. Then, periodic points for $\mathcal P_\xi$ of principal period $n$ correspond to functions with exactly $n$ local extrema in one period.

Another property of this section is that $\mathrm{Fix}(R)\subset \Pi$, where 
\begin{equation*}
	R(x,y,z,w) = (x,-y,z,-w).
\end{equation*}
is a reversing symmetry of \eqref{eq:SH_system}. It is well known \cite{lamb1992reversing, wilczak2003chaos} that in such case the Poincar\'e map is also $R$-reversible. In order to find an even ($R$-symmetric) periodic orbit of \eqref{eq:SH_system} it suffices to find a point $u=(x,0,z,0)\in\mathrm{Fix}(R)$ such that $\mathcal P_\xi^k(u)\in\mathrm{Fix}(R)$. Then $\mathcal P_\xi^{2k}(u)=u$ and the trajectory of $u$ is periodic and $R$-symmetric.

Hence, even periodic solutions of \eqref{eq:SH_system} correspond to solutions of the scalar equation
\begin{equation}\label{eq:tildeG}
	\pi_w \mathcal P^2_\xi(x,0,z,0) = 0,
\end{equation}
where $\pi_w$ is the projection onto $w$ variable. Observe that the intersection of Poincar\'e section with the energy level $E=0$
$$E(x,0,z,w) = -\frac{1}{2}z^2+\frac{1}{4}(x^2-1)=0$$ 
gives the following relation
\begin{equation}\label{eq:Z}
	u''=z = \pm \frac{1}{\sqrt{2}}(x^2 - 1).
\end{equation}
This means, that a solution $u$ at the energy level $E=0$ has always proper extrema provided they are taken at $x\neq\pm1$. Shifting minumum $\tilde u_1<-1$ to $t=0$ (see geometric conditions in Theorem~\ref{sec:forcingTheorem} and Fig.~\ref{fig:periodic_orbits}) we must take $z(0)=\frac{1}{\sqrt2}(x(0)^2-1)>0$. Substituting \eqref{eq:Z} to \eqref{eq:tildeG}, we eventually transform the question of finding even periodic solutions of \eqref{eq:SH_system} to a zero finding problem of the following equation
\begin{equation}\label{eq:equationForPOs}
	G(\xi,x) = \left (\pi_w \circ \mathcal P_\xi^2 \right ) \left (x,0, \frac{1}{\sqrt{2}}(x^2 - 1), 0 \right) = 0.
\end{equation}

Finally, we have to check that an orbit corresponding to a solution of $G(\xi,x_0)=0$ satisfies the geometric conditions from Theorem~\ref{sec:forcingTheorem}. It is sufficient to check if
\begin{equation}\label{eq:extremaConstraint}
	(x_0 < -1) \wedge (x_1 > 1)\wedge (-1<x_2<1),
\end{equation}
where 
$x_i= \pi_x\mathcal P_\xi^i(x_0,0,\frac{1}{\sqrt{2}}(x_0^2 - 1),0)$, $i=1,2$. Indeed, by \eqref{eq:equationForPOs} the function has at most four extrema in one period. From \eqref{eq:extremaConstraint} the minima $x_0\neq x_2$ are different. Hence, the solution has exactly four extrema in period. All of them are proper because $x_i\neq\pm1$, $i=0,1,2$.

\section{Computer-assisted proof of Theorem~\ref{thm:main}}\label{sec:cap}
\begin{figure}[ht] 
	\centerline{\includegraphics[width=.95\textwidth]{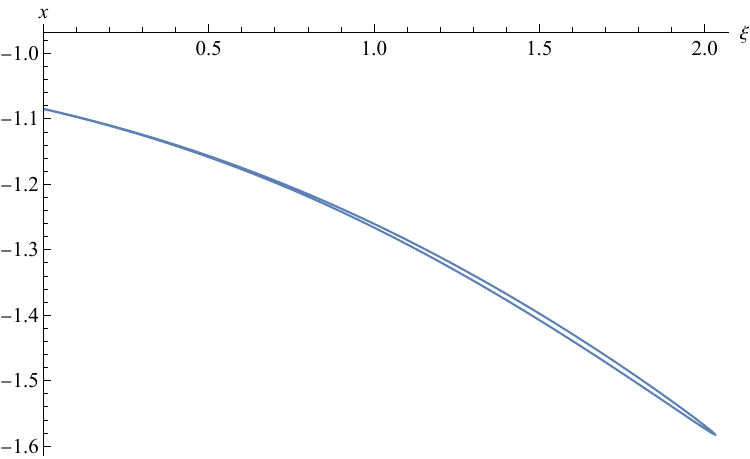}}
	\caption{Branch of even periodic orbits on the energy level $E=0$. At $\xi^*\approx 2.0316516135713902$ the fold bifurcation occurs.}
	\label{fig:bifdiag}
\end{figure}

\begin{figure}[ht] 
	\centerline{\includegraphics[width=.95\textwidth]{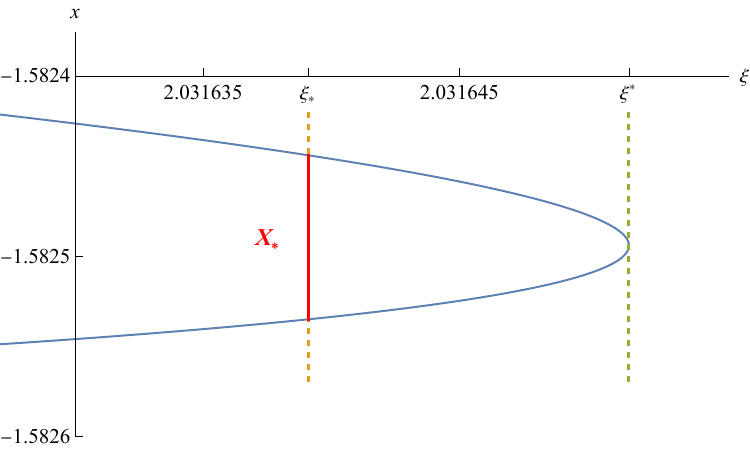}}
	\caption{Branch of even periodic orbits on the energy level $E=0$. At $\xi^*\approx 2.0316516135713902$ the fold bifurcation occurs. $\xi_*$ is a threshold value at which we switch parameterisation of this curve. For $\xi\geq \xit$ the curve is parameterised as $\tilde \xi:X_*\to\R$ and the domain $X_*$ contains both values $x_\pm(\xit)$.}
	\label{fig:fold}
\end{figure}

As already mentioned (see \eqref{eq:changeVariables}) the parameter range $\alpha\geq\frac{3}{2}$ corresponds to $0<\xi\leq \sqrt 8$. Numerical simulation shows (see Fig.~\ref{fig:bifdiag}, Fig.~\ref{fig:fold} and \cite{BurkeKnoblach}) that a saddle-node bifurcation occurs at $\xi\approx 2.0316516135713902$ creating two branches of even periodic orbits. These branches continue to exist until $\xi=0$ (in fact they continue until $\xi\approx -2.3$ but we restrict here to the range $\xi \geq 0$ due to change of variables \eqref{eq:changeVariables}).

Computer-assisted proof of Theorem~\ref{thm:main} is split into the following steps.
\begin{description}
	\item[Step 1.] We fix a threshold value $\xit=266\,291 \cdot 2^{-17} \approx2.0316390991210938$ ($\xit$ is IEEE-754 \cite{IEEE754-2019} representable number) and prove that there are two smooth branches $x_\pm:[0,\xit]\to\mathbb R$ such that for all $\xi\in[0,\xit]$ there holds $x_-(\xi)<x_+(\xi)$, $G(\xi,x_\pm(\xi))\equiv 0$ and $x_{\pm}(\xi)$ satisfies \eqref{eq:extremaConstraint}.
	\item [Step 2.] We prove that there is a smooth and concave function $\tilde\xi:[x_-(\xit),x_+(\xit)]\to\mathbb R$ such that for $x\in [x_-(\xit),x_+(\xit)]$ there holds $G(\tilde \xi(x),x)\equiv 0$ and $x$ satisfies \eqref{eq:extremaConstraint} in the system with parameter $\tilde \xi(x)$.
	\item [Step 3.] We prove that $\tilde \xi$ has unique maximum $\xi^*=\tilde \xi(x^*)$ for some $x^*\in[x_-(\xit),x_+(\xit)]$. 
	\item [Step 4.] We prove that graphs of $x_\pm$ and $\tilde \xi$ glue into a smooth curve in the $(\xi,x)$ plane, as shown in Fig.~\ref{fig:bifdiag} and Fig.~\ref{fig:fold}.
\end{description}

Finally, from these steps we will conclude that the functions $x_\pm$ can be extended beyond $\xit$ to $\xi^*$, at which value they are equal and the saddle-node bifurcation occurs at $(\xi^*,x^*)$. Via change of variables \eqref{eq:changeVariables} we obtain the assertion of Theorem~\ref{thm:main}.

\subsection{Interval Newton Operator}
In a computer-assisted verification of \textbf{Steps 1--4} we will use the standard Interval Newton Operator for proving the existence and uniqueness of zeros of maps \cite{Neumaier1990}.
\begin{theorem}\label{thm:INO}
	Let $f:X\subset \mathbb R^n\to\mathbb R^n$ be $\mathcal C^1$ smooth, where $X$ is compact and convex. Fix $x_0\in\mathrm{int} X$ and let $[A]$ be an interval matrix such that 
	\begin{eqnarray*}
		\{Df(x) : x\in X\} &\subset& [A].
	\end{eqnarray*}
	If the interval Newton operator
	$$
	N := x_0 - [A]^{-1}f(x_0)
	$$
	is well defined and $N\subset \mathrm{int}X$ then $f$ has unique zero $\widetilde x$ in $X$ and $\widetilde x\in N$.
\end{theorem}
Using implicit function theorem and Theorem~\ref{thm:INO} it is easy to prove the following extension of Theorem~\ref{thm:INO} to parameterised functions.

\begin{theorem}\label{thm:parameterisedINO}
	Let $f:Z\times X\subset \mathbb R^k\times\mathbb R^n\to\mathbb R^n$ be $\mathcal C^1$ smooth, where $Z,X$ are compact and convex. Fix $x_0\in\mathrm{int} X$ and let $[A]$ be an interval matrix and $[e]$ be an interval vector such that 
	\begin{eqnarray*}
		\{D_xf(z,x) : (z,x)\in Z\times X\} &\subset& [A],\\
		\{f(z,x_0) : z\in Z\} &\subset & [e].
	\end{eqnarray*}
	If the interval Newton operator
	$$
	N := x_0 - [A]^{-1}[e]
	$$
	is well defined and $N\subset \mathrm{int}X$ then the set of zeroes of $f$ in $Z\times X$ is a graph of a $\mathcal C^1$ smooth function $\widetilde x:Z\to N\subset X$. That is, for $(z,x)\in Z\times X$
	$$
	f(z,x)=0\Longleftrightarrow x=\widetilde x(z).
	$$
\end{theorem}

\subsection{Details of verification of Steps 1--4.}

\textbf{Step 1.} To prove the existence of two smooth curves  $x_\pm:[0,\xit]\to\mathbb R$ solving \eqref{eq:equationForPOs} we perform an adaptive subdivision of the parameter range $[0,\xit]=\Xi^1_\pm\cup\cdots\cup \Xi_{\pm}^{N_\pm}$, where $\Xi^j_\pm$ are closed intervals satisfying $\min\Xi_\pm^1=0$, $\max \Xi_\pm^{N_\pm}=\xit$ and $\max \Xi_\pm^j=\min \Xi_\pm^{j+1}$ for $j=1,\ldots,N_{\pm}-1$. The sizes of subdivisions $N_-=24157$, $N_+=142821$ are returned by our adaptive algorithm. Then for each subinterval $\Xi^j_\pm$ of parameters we construct (using continuation algorithms) a closed interval $X_{\pm}^j$ and using Theorem~\ref{thm:parameterisedINO} we prove that the solution set of \eqref{eq:equationForPOs} in $\Xi^j_\pm\times X_{\pm}^j$ is a smooth curve $x_\pm^j:\Xi_\pm^j\to X_\pm^j$. Moreover \eqref{eq:extremaConstraint} holds true for every point on this curve. Finally we show, that 
\begin{eqnarray*}
	X_\pm^j\cap X_\pm^{j+1}&\neq& \emptyset,\quad j=1\ldots,N_\pm-1.
\end{eqnarray*}
Because $\max \Xi^j_\pm=\min\Xi_\pm^{j+1}$ and by the uniqueness property of the Interval Newton Operator we conclude that $x_\pm^j\left(\max\Xi_{\pm}^j\right)=x_\pm^{j+1}\left(\min\Xi_{\pm}^{j+1}\right)$ for $j=1,\ldots,N_\pm-1$. Hence, the functions $x_\pm^j$, defined on subintervals $\Xi_\pm^j$ of the parameter range join into continuous functions $x_\pm$ defined on $[0,\xit]$. Sine at every $\xi\in[0,\xit]$ we have $D_xG(\xi,x_\pm(\xi))\neq 0$ (the interval Newton operator is defined), by the implicit function theorem these functions are smooth. 

Finally, using Theorem~\ref{thm:INO} we validate that
\begin{equation}\label{eq:xitBound}
	x_+(\xit)\in -1.5824440318_{2912}^{327613},\quad     
	x_-(\xit)\in -1.58253506275_{30319}^{63035}.
\end{equation}
From these bounds it is clear, that $x_-(\xit)< x_+(\xit)$ and thus, by the uniqueness property of the interval Newton operator this relation extends to $x_-(\xi) < x_+(\xi)$ for all $\xi\in[0,\xit]$.

\textbf{Step 2.}
To prove the existence of a concave curve $\tilde{\xi} : [x_{-}(\xit), x_{+}(\xit)] \rightarrow \R$ solving \eqref{eq:equationForPOs} we need to compute second order derivatives of $G$ with respect to the parameter $\xi$. Because of the interface of the CAPD library \cite{capd}, in order to compute derivatives with respect to a parameter we have to treat it as a variable. Therefore we extend the system \eqref{eq:SH_system} to
\begin{equation} \label{eq:SH_extended_system}
	\begin{cases}
		x'(t) &= y(t) \\
		y'(t) &= z(t) \\
		z'(t) &= w(t) \\
		w'(t) &= -\xi(t) z(t) + x(t) - x^3(t) \\
		\xi'(t) &= 0.
	\end{cases}
\end{equation}
We are going to seek the zeros of the modified function
\begin{equation} \label{eq:equationForPOs_extended}
	\bar{G}(\xi,x) = \left (\pi_w \circ \mathcal P^2 \right ) \left (x,0, \frac{1}{\sqrt{2}}(x^2 - 1), 0, \xi \right),
\end{equation}
where $\mathcal{P} : \bar{\Pi} \rightarrow \bar{\Pi}$ is the Poincar\'e map for the section
\begin{equation*}
	\bar{\Pi} = \left \{(x,y,z,w,\xi) : y = 0 \right \}.
\end{equation*}
Clearly zeroes of $G$ are in on-to-one correspondence with zeroes of $\bar G$.

The rest of this step is similar to the \textbf{Step 1} with the addition of verification of concavity. Denote by $x^{\xit}_{-}$ the lower bound for $x_{-}(\xit)$ and by $x^{\xit}_{+}$ the upper bound for $x_{+}(\xit)$ -- see \eqref{eq:xitBound}. To find a curve solving \eqref{eq:equationForPOs_extended} we perform an adaptive subdivision of the range $X_*:=[x^{\xit}_{-}, x^{\xit}_{+}] = \bigcup_{j=1}^{N} X^j$ into $N = 6903$ closed intervals with $\max X^j=\min X^{j+1}$, $j=1,\ldots, N-1$. Then for each subinterval $X^j$ we construct an interval $\Xi^{j}$ and using Theorem~\ref{thm:parameterisedINO} we prove, that the solution set of \eqref{eq:equationForPOs_extended} in $\Xi^{j} \times X^{j}$ is the graph of a smooth function $\tilde \xi^j:\Xi^{j} \to X^{j}$. We also check $\Xi^{j} \cap \Xi^{j + 1}\neq\emptyset$ for $j = 1,...,N - 1$. By the uniqueness property of the interval Newton operator the functions $\tilde \xi^j$ join into a smooth function $\tilde \xi:X_*\to\R$ such that for $x\in X_*$ there holds $\bar G(\tilde \xi(x),x)=0$ and $x$ satisfies \eqref{eq:extremaConstraint}. 

Additionally, on each subinterval $X^j$ we compute a bound on the second order derivative of the implicit function $\tilde \xi$. For this purpose we use the $\mathcal C^r$-Lohner algorithm \cite{WilczakZgliczynski2011} for integration of second order variational equations for \eqref{eq:SH_extended_system} needed for the second derivatives of Poincar\'e map $\mathcal P$. From these computation we obtain a bound
\begin{equation*}
	\tilde{\xi}''(x) \in  [-74010.849232287583, -12744.872650106316],
\end{equation*}
for $x\in X_*$. Hence the function $\tilde{\xi}$ is concave on the interval $[x^{\xit}_{-},x^{\xit}_{+}]$.

\textbf{Step 3.} We apply Theorem~\ref{thm:INO} to the following equation
$$
H(\xi,x) = \left(G(\xi,x),\frac{\partial G}{\partial x}(\xi,x)\right)=0.
$$
Using standard Newton method we have found an approximate solution
$$(\xi_0,x_0) = (2.0316516135713902,-1.5824941113082425)$$
to $H(\xi,x)=0$. Then, using Theorem~\ref{thm:INO} we proved that there is a unique zero $(\xi^*,x^*)$ of $H$ in the set $(\xi_0,x_0)+[-1,1]^2\cdot 10^{-10}$ that belongs to 
\begin{eqnarray}
	\xi^*&\in& 2.0316516135_{613893}^{814116}, \label{eq:xiBifBound}\\
	x^*&\in& -1.582494111^{3301776}_{2863635}.\label{eq:xBifBound}
\end{eqnarray}
Using bound \eqref{eq:xiBifBound} and change of variables \eqref{eq:changeVariables} we obtain a bound for $\alpha^*$ in Theorem~\ref{thm:main}. Finally using bounds \eqref{eq:xiBifBound}-\eqref{eq:xBifBound} we checked, that $(\xi^*,x^*)$ belongs to one of the sets $\Xi^j\times X^j$ used to verify the existence of $\tilde \xi$ in \textbf{Step 2}. Hence, by the uniqueness property of the interval Newton operator, we obtain $\xi^*=\tilde \xi(x^*)$. From the implicit function theorem we have $\tilde\xi'(x^*)=-\frac{\partial G}{\partial x}(\xi^*,x^*)/\frac{\partial G}{\partial \xi}(\xi^*,x^*)=0$. Given that $\tilde \xi$ is concave we conclude, that $\tilde \xi$ has unique maximum at $\xi^*$.

\textbf{Step 4.}
By the construction of $\tilde \xi$ we know that $x_\pm(\xit)\in X_*$.
To show, that the graphs $x_{\pm}$ and $\tilde{\xi}$ glue into a smooth curve in the $(\xi,x)$ plane it suffices to check that 
\begin{equation}\label{eq:glueCondition}
	\tilde\xi(x_\pm(\xit))=\xit.
\end{equation}
Indeed, because the interval Newton operator is well defined in \textbf{Steps 1--2}, both partial derivatives $\frac{\partial G}{\partial x}(\xit,x_\pm(\xit))$ and $ \frac{\partial G}{\partial \xi}(\xit,x_\pm(\xit))$ are nonzero. Hence the curves $x_{\pm}(\xit)$ can be extended in a smooth and unique way in a neighborhood of $\xit$ and therefore graphs of $\tilde \xi$ and $x_\pm$ must locally coincide near $(\xit,x_\pm(\xit))$.

In order to check \eqref{eq:glueCondition}, using bounds \eqref{eq:xitBound} we verify that 
$$(\xit,x_-(\xit))\in \Xi^1\times X^1\quad\text{and}\quad (\xit,x_+(\xit))\in \Xi^N\times X^N.$$
Then, from construction of $\tilde \xi$ in \textbf{Step 2} we obtain \eqref{eq:glueCondition}.

\subsection{Implementation notes.}
In the computer-assisted verification of \textbf{Steps 1--4} we need to compute bounds on Poincar\'e map and its derivatives up to the second order. For this purpose we used interval arithmetics \cite{Moore1966,IEEE1788-2015}, algorithms for rigorous integration of ODEs and associated (higher order) variational equations and algorithms computation of Poincar\'e maps \cite{poincare} all implemented in the CAPD library \cite{capd}.

The C++ program that performs verification of \textbf{Steps 1--4} is a supplement to this article and available at \cite{repo}. 

The computation of the bound for the curve $x_{-}(\xi)$ takes approximately 109 seconds. This curve corresponds to the solution from \cite{chaos}. The curve $x_{+}(\xi)$ takes less than 10 minutes (537 seconds) to compute, and the curve $\tilde{\xi}(x)$ 1.5 minutes (100 seconds). Overall, the total computation takes less than 13 minutes.

\section{Conclusions and future works}
In this paper we reproduced and extended some results from \cite{chaos} about dynamics of the stationary Swift-Hohenberg equation on the energy level $E=0$. After change of variables \eqref{eq:changeVariables} we obtain a system with a parameter $\xi$ related to original parameter $\alpha$. Numerical simulation shows that the two curves resulting from Theorem~\ref{thm:main} continue to exist for negative values of $\xi$, which apparently undergo further bifurcation near $\xi\approx -2.3$. 

Our method is computationally less expensive than the one proposed in \cite{chaos}. This gives a hope to obtain results as in Theorem~\ref{thm:main} for a range of energy levels and, perhaps, study codimension two bifurcations of this family of periodic orbits.

\end{document}